\documentclass{amsart}
\usepackage{amssymb,latexsym,amsmath,amsfonts,hhline,color}
\usepackage{pdfsync}

\newcommand{\cal}{\mathcal}

\makeatletter
\renewcommand{\subsection}{\@startsection{subsection}{2}{0mm}{-2mm}{-2mm}{\bf\normalsize}}

\def\sbsnt#1{\subsection{#1}}
\makeatother

\newtheorem{formula}{}[section]
\newtheorem{definition}[formula]{Definition}
\newtheorem{corollary}[formula]{Corollary}
\newtheorem{remark}[formula]{Remark}
\newtheorem{proposition}[formula]{Proposition}
\newtheorem{lemma}[formula]{Lemma}
\newtheorem{theorem}[formula]{Theorem}

\def\thrm{\begin{theorem}}
\def\thrml#1{\begin{theorem}\label{#1}}
\def\ethrm{\end{theorem}}
\def\rmrk{\begin{remark}}
\def\rmrkl#1{\begin{remark}\label{#1}}
\def\ermrk{\end{remark}}
\def\prpstn{\begin{proposition}}
\def\prpstn#1{\begin{proposition}\label{#1}}
\def\eprpstn{\end{proposition}}
\def\dfntn{\begin{definition}}
\def\dfntnl#1{\begin{definition}\label{#1}}
\def\edfntn{\end{definition}}
\def\nmrt{\begin{enumerate}}
\def\enmrt{\end{enumerate}}
\def\tm#1{\item[{\rm (#1)}]}
\def\qtn{\begin{equation}}
\def\qtnl#1{\begin{equation}\label{#1}}
\def\eqtn{\end{equation}}
\def\lmm{\begin{lemma}}
\def\lmml#1{\begin{lemma}\label{#1}}
\def\elmm{\end{lemma}}
\def\crllr{\begin{corollary}}
\def\crllrl#1{\begin{corollary}\label{#1}}
\def\ecrllr{\end{corollary}}
\def\css{\begin{cases}}
\def\ecss{\end{cases}}

\def\proof{\noindent{\bf Proof}.\ }

\def\cA{{\cal A}}

\def\cS{{\cal S}}

\def\cX{{\cal X}}

\def\cZ{{\cal Z}}

\def\mC{{\mathbb C}}

\def\mQ{{\mathbb Q}}
\def\mZ{{\mathbb Z}}

\DeclareMathOperator{\aut}{Aut}

\DeclareMathOperator{\cla}{Cla}

\DeclareMathOperator{\Inn}{Inn}

\DeclareMathOperator{\Irr}{Irr}

\DeclareMathOperator{\ORT}{O}

\DeclareMathOperator{\rad}{rad}
\DeclareMathOperator{\rk}{rk}

\DeclareMathOperator{\Span}{Span}
\DeclareMathOperator{\SP}{Sp}

\DeclareMathOperator{\sym}{Sym}
\DeclareMathOperator{\tr}{tr}

\def\bul{\hfill\vrule height .9ex width .8ex depth -.1ex}
\def\bull{\hfill\vrule height .9ex width .8ex depth -.1ex\medskip}

\def\qaq{\quad\text{and}\quad}
\def\qoq{\quad\text{or}\quad}
\def\mmod#1#2#3{#1=#2\ (\text{\rm mod}\hspace{2pt}#3)}

\newcommand{\grp}[1]{\langle {#1}\rangle}
\newcommand{\und}[1]{{\underline{#1}}}

\begin{document}
\title{The Schur-Wielandt theory for central S-rings}
\author{Gang Chen}
\address{School of Mathematics and Statistics, Central China Normal University, Wuhan, China}
\email{chengang19762002@aliyun.com}
\author{Mikhail Muzychuk}
\address{Netanya Academic College, Netanya, Israel}
\email{muzy@netanya.ac.il}
\author{Ilya Ponomarenko}
\address{Steklov Institute of Mathematics at St. Petersburg, Russia}
\email{inp@pdmi.ras.ru}
\thanks{The work of the first author was Financially supported by self-determined research funds of CCNU (No.CCNU15A02031) from the colleges’basic research and operation of MOE. The work of the third author was partially supported by the RFBR Grant 14-01-00156}
\date{}

\maketitle

\begin{abstract}
Two basic results on the S-rings over an abelian group are the Schur theorem on multipliers and the Wielandt
theorem on primitive S-rings over groups with a cyclic Sylow subgroup. None of these theorems is directly generalized
to the non-abelian case. Nevertheless, we prove that they are true for the central S-rings, i.e., for those which are
contained in the center of the group ring of the underlying group (such S-rings naturally arise in the supercharacter
theory). We also generalize the concept of a B-group introduced by Wielandt, and show that any Camina group is
a generalized B-group whereas with few exceptions, no simple group is of this type.
\end{abstract}

\section{Introduction}\label{sec:1}

A {\it Schur ring} or {\it S-ring} over a finite group $G$ can be defined as a subring of the group ring $\mZ G$
that is a free $\mZ$-module spanned by a partition of $G$ closed under taking inverse and containing the identity $e$
of $G$ as a class (see Section~\ref{150315x} for details). The S-ring theory was initiated by Schur~\cite{S} and then
developed by Wielandt~\cite{Wie} who wrote in \cite{Wie69}  that S-rings provide one ``of three major tools'' to study
a group action.\footnote{The two other tools are the representation theory and the method of invariant relations.}\medskip 

Until recently, the focus was on studying S-rings over abelian groups and the main applications of this theory were
connected with algebraic combinatorics problems~\cite{MP}. However, as it was observed  in \cite{He10}, the 
supercharacter theory developed to study group representations, is nothing else than the theory of commutative
S-rings of a special form that we call here {\it central}.

\dfntn
An S-ring over a group $G$ is said to be {\it central} if  it is contained in the center $\cZ(\mZ G)$ 
of the group ring $\mZ G$.  
\edfntn

 An example of such a ring is obtained from any permutation group $K$ such that
$$
G\Inn(G)\le K\le \sym(G),
$$
where $\Inn(G)$ is  the inner automorphism group of $G$; the corresponding 
partition of~$G$ is formed by the orbits of the  stabilizer of $e$ in~$K$. In the special case when $K=\sym(G)$,
this produces the {\it trivial} central S-ring $\mZ e+\mZ\und{G}$, where $\und{G}$ is the sum of all elements of $G$. 
On the other hand, if  $K=G\Inn(G)$, the orbits are 
the conjugacy classes of $G$; this shows that $\cZ(\mZ G)$ is a central S-ring.
In particular, any S-ring over an abelian group is central. The main goal of the present paper is to extend the basic results 
on S-rings from abelian case to the central one.\medskip

The Schur theorem on multipliers is a fundamental statement in the theory of S-rings over abelian groups. To explain it,
given an integer $m$ coprime to $|G|$, we define a permutation on the elements of the group $G$ by
$$
\sigma_m:G\to G,\ x\mapsto x^m.
$$  
It  permutes also the conjugacy classes of $G$, and so induces a linear isomorphism of  the ring $\cZ(\mZ G)$.  If 
the group~$G$ is abelian, then $\sigma_m\in\aut(G)$, $\cZ(\mZ G)=\mZ G$ and the Schur theorem on multipliers states
that $\sigma_m$ is a Cayley automorphism of every S-ring over $G$. Our first result shows that in the nonabelian case,
$\sigma_m$ is still an automorphism (but not a Cayley one) of  any central S-ring over~$G$.

\thrml{100315v}
Let $\cA$ be a central S-ring over a group $G$, and let $m$ be an integer coprime to $|G|$. Then $\sigma_m(\cA)=\cA$ and 
$\sigma_m|_\cA\in\aut(\cA)$.
\ethrm

Based on this result  for the abelian case, Wielandt generalized the Schur theorem on primitive groups having a regular cyclic subgroup. In fact,
the Wielandt proof shows that if $G$ is an abelian group of composite order that has a cyclic Sylow subgroup, then
no proper S-ring over $G$ is primitive.\footnote{The primitivity concept in S-ring theory plays the same role as the simplicity in group theory.} 
The following statement establishes ``a central version'' of the Wielandt theorem.

\thrml{070215a}
Let $\cA$ be a nontrivial central S-ring over a group $G$ of composite order. Suppose that $G$ has a normal cyclic 
Sylow $p$-subgroup. Then $\cA$ is imprimitive.
\ethrm

 Following \cite{Wie}, a finite group $G$ is called a B-group if every primitive group containing a regular subgroup 
isomorphic to~$G$ is 2-transitive.  It should be remarked that most of the B-groups $G$ mentioned in~\cite{Wie} 
satisfy a priori a stronger condition: no nontrivial S-ring over $G$ is primitive. In this sense, the following definition seems to be quite natural.
In what follows, we say that a central S-ring over $G$ is {\it proper} if it lies strictly between  $\cZ(\mZ G)$ and the trivial S-ring over~$G$.

\dfntn
A group $G$ is called a generalized B-group if no proper central S-ring over $G$ is primitive. 
\edfntn

Clearly, every B-group is also a generalized one. The converse statement is not true; see Subsection~\ref{190315a}.
A nontrivial example of a generalized B-group is given in Theorem~\ref{070215a}.
The following statement gives a family of  generalized B-groups; we don't know whether they 
are B-groups. Below, under a {\it Camina} group,  we mean a group $G$  that has a proper 
nontrivial  normal subgroup~$H$  such that  each  $H$-coset    distinct  from $H$  is  contained in a conjugacy 
class of~$G$ (in other terms, $(G,H)$ is a Camina pair).\footnote{To simplify the presentation, we use the term ``Camina group" not only in the case
where $(G,G')$ is a Camina pair.}

\thrml{100315a}
Any  Camina group is  a generalized B-group. 
\ethrm

The class of the Camina groups includes, in particular, all Frobenius  and extra-special groups; see~\cite{Ca}. 
Thus, by Theorem~\ref{100315a}, we obtain the following statement.

\crllrl{100315u}
Any Frobenius or extra-special group  is  a generalized B-group.\bull
\ecrllr
 
The last result of the present paper shows that with  a few possible exceptions, no simple 
group is a generalized B-group. The proof  is based on the Schur theorem on multipliers and 
the characterization of rational simple groups given in~\cite{FS}.

\thrml{070215b}
A generalized B-group $G$ is not simple  unless $|G|\le 3$, or $G\cong\SP(6,2)$ 
or $\ORT^+(8,2)'$.\footnote{In fact, we do not know whether two simple groups from Theorem~\ref{070215b} are generalized B-groups.}
\ethrm

For the reader convenience, we collect the basic facts on S-rings in Section~\ref{150315x}. The proofs of 
Theorems~\ref{100315v} and~\ref{070215a} are contained in Sections~\ref{190315t} and \ref {190315u}, respectively. 
The results  concerning generalized B-groups are in Section~\ref{190315v}.\medskip

{\bf Notation.}

As usual, $\mZ$, $\mQ$ and $\mC$ denote the ring of integers and the fields of rationals and complex numbers, respectively.

The identity of a group $G$ is denoted by $e$; the set of non-identity elements in $G$ is denoted by  $G^\#$.

The set of conjugacy classes of $G$ is denoted by $\cla(G)$.

Let $X\subseteq G$. The subgroup of $G$ generated by $X$ is denoted by $\grp{X}$; 
we also set $\rad(X)=\{g\in G:\ gX=Xg=X\}$. 

The element $\sum_{x\in X}x$ of the group  ring $\mZ G$ is denoted by $\und{X}$. 

For an integer $m$, we set $X^{(m)}=\{x^m:\ x\in X\}$ and $\und{X}^{(m)}=\und{X^{(m)}}$.

The group of all permutations of the elements of $G$ is denoted by $\sym(G)$. 

The additive and multiplicative groups of the ring $\mZ/(n)$ are denoted by  $\mZ_n$ and $\mZ^*_n$, respectively.

\section{Preliminaries}\label{150315x}

Let $G$ be a finite group. A subring~$\cA$ of the group ring~$\mZ G$ is called a {\it Schur 
ring} ({\it S-ring}, for short) over~$G$ if there exists a partition $\cS=\cS(\cA)$ of~$G$ 
such that
\nmrt
\tm{S1} $\{e\}\in\cS$,
\tm{S2} $X\in\cS\ \Rightarrow\ X^{-1}\in\cS$,
\tm{S3} $\cA=\Span\{\und{X}:\ X\in\cS\}$.
\enmrt
In particular, for all $X,Y,Z\in\cS$ there is a nonnegative integer $c_{XY}^Z$ such that
$$
\und{X}\,\und{Y}=\sum_{Z\in\cS}c_{XY}^Z\und{Z},
$$
these integers are the structure constants of $\cA$ with respect to the linear base $\{\und{X}:\ X\in\cS\}$.
The number $\rk(\cA)=|\cS|$ is called the {\it rank} of~$\cA$.\medskip

Let $\cA'$ be an S-ring over a group $G'$. Under a {\it Cayley isomorphism} from $\cA$ to~$\cA'$, we mean
a group isomorphism $f:G\to G'$ such that $\cS(\cA)^f=\cS(\cA')$. This is a special case of the ordinary {\it isomorphism};
by definition, it is a  bijection $f:G\to G'$ that induces a ring isomorphism from $\cA$ to $\cA'$ taking
$\und{X}$ to $\und{X'}$ for all $X\in\cS$,  where $X'=X^f$.\medskip

The classes of the partition $\cS$ are called the {\it basic sets} of the S-ring~$\cA$. Any union of them
is called an {\it $\cA$-set}. Thus, $X\subseteq G$ is an $\cA$-set if and only if $\und{X}\in\cA$. The set of all $\cA$-sets
is closed with respect to taking inverse and product. Any subgroup of~$G$ that is an $\cA$-set, is called an {\it $\cA$-subgroup} 
of~$G$ or {\it $\cA$-group}. With each $\cA$-set $X$, one can  naturally associate two $\cA$-groups, 
namely $\grp{X}$ and  $\rad(X)$ (see Notation). The S-ring $\cA$ is called {\it primitive} if the only $\cA$-groups are $e$ and $G$,
otherwise this ring is called  {\it imprimitive}.\medskip

We will use the following statement proved in \cite[Proposition~22.3]{Wie}. Below for a function $f:\mZ\to \mZ$ and an
element $\xi=\sum_ga_gg$ of the ring $\mZ G$, we set $f[\xi]=\sum_gf(a_g)g$.

\lmml{110315c}
Let $\cA$ be an S-ring, $f:\mZ\to \mZ$ an arbitrary function and $\xi\in\cA$. Then $f[\xi]\in\cA$.\bull 
\elmm

The important special case is  when $f(a)=1$ or $0$ depending on whether $a\ne 0$ or $a=0$. Then
$f[\xi]=\und{X}$, where $X$ is the support of~$\xi$, and we refer to Lemma~\ref{110315c}
as to the Schur-Wielandt principle.

\section{The Schur theorem on multipliers}\label{190315t}

{\bf Proof of Theorem~\ref{100315v}.} Since obviously $\sigma_{mm'}=\sigma_{m^{}}\sigma_{m'}$ for
all $m$ and $m'$, without loss of generality, we can assume that $m$ is a prime. We need the following auxiliary 
lemma.

\lmml{100315b}
Let $X\in\cla(G)$ and $p$ an arbitrary prime. Then 
$$
\und{X}^p=\sum_{Y\in\cla(G)}a_Y\und{Y}
$$ 
for some nonnegative integers $a_Y$'s. Moreover,
$$
a_Y|Y|=\css
|X|\ (\text{\rm mod}\hspace{2pt}p), &\text{if $Y=X^{(p)}$,}\\
\ 0\quad (\text{\rm mod}\hspace{2pt}p), &\text{if $Y\ne X^{(p)}$.}\\
\ecss
$$
\elmm
\proof The first statement follows from the fact that $\cZ(\mZ G)$ is an S-ring.  To prove the second one, set
$$
T_Y=\{(x_1,\ldots,x_p)\in X^p:\ x_1\cdots x_p\in Y\}.
$$
Clearly, $(T_Y)^G=T_Y$. Since also $(x_1x_2\cdots x_p)^{x_1}=x_2\cdots x_px_1$, the set $T_Y$ is invariant
with respect to the cyclic shift  $\pi:(x_1,x_2,\ldots,x_p)\mapsto (x_2,\ldots,x_p,x_1)$. Moreover, 
since $p$ is prime, we have 
$$
|(x_1,\ldots,x_p)^{\grp{\pi}}|=1\ \text{or}\ p.
$$
However,  $|(x_1,\ldots,x_p)^{\grp{\pi}}|=1$  if and only if $x_1=\cdots=x_p$, which is possible only if $Y=X^{(p)}$;
in the latter case, the group $\grp{\pi}$ has exactly $|X|$ orbits of the form $\{(x,...,x)\}$, $x\in X$. Taking into
account that $T_Y$ is a disjoint union of $\grp{\pi}$-orbits, we have
$$
|T_Y|=|X|\delta+pu
$$
where $\delta=\delta_{Y,X^{(p)}}$ is the Kronecker delta and $u$ is the number of the $\grp{\pi}$-orbits of size $p$.
Thus, the required statement follows because $a_Y|Y|=|T_Y|$.\bull 

Let us continue the proof of the theorem. Since $p:=m$ is coprime to $|G|$, the mapping 
$$
\und{X}\mapsto\und{X}^{(p)},\quad X\in\cla(G),
$$ 
is a bijection. It  
induces a linear isomorphism of the ring $\cZ(\mZ G)$; the image of the element $\xi$ under this isomorphism is 
denoted by $\xi^{(p)}$. From Lemma~\ref{100315b}, it follows that $\mmod{\und{X}^p}{\und{X}^{(p)}}{p}$; here,
we make use of the fact that $|X^{(p)}|=|X|$ for all~$X$. Therefore,
\qtnl{110315a}
\xi^{(p)}=f[\xi^p]\quad\text{for all}\quad\xi\in\cZ(\mZ G).
\eqtn
where $f(a)$ is the remainder in the division of $a$ by $p$.\medskip

To prove the first part of the theorem, let $X\in\cS(\cA)$. Then $X$ is a union of some classes $X_i\in\cla(G)$, $i\in I$. Thus, 
by~\eqref{110315a}, we have
\qtnl{120315c}
\sigma_p(\und{X})=\und{X}^{(p)}=\sum_{i\in I}\und{X_i}^{(p)}=\sum_{i\in I}f[\und{X_i}^p]=f[\sum_{i\in I}\und{X_i}^p]=
f[(\sum_{i\in I}\und{X_i})^p]=f[\und{X}^p].
\eqtn
However, by Lemma~\ref{110315c}, the right-hand side belongs to~$\cA$. Therefore $\sigma_p(\und{X})\in\cA$ and $X^{(p)}$
is an $\cA$-set. Moreover, suppose that  it contains a proper basic set $Y$. By the Dirichlet Theorem, one can find a
prime  $p'$  such that $\mmod{pp'}{1}{n}$, where $n=|G|$.  Now, the  above argument shows that
$Y^{(p')}$ is a proper $\cA$-subset of $X$. Thus $X^{(p)}\in\cS(\cA)$ and so $\sigma_p(\cA)=\cA$.\medskip

To prove the second part of the theorem, it suffices to verify that $\sigma_m$ induces a ring isomorphism of  $\cZ(\mZ G)$:
then it is, obviously, an S-ring isomorphism of $\cZ(\mZ G)$ that takes $\cA$ to itself, and hence it is an isomorphism
of $\cA$, as required. To do this, without loss of generality, we can assume that $p>2n$ (for otherwise,
by the Legendre theorem, there exists  a prime $q>2n$ such that $\mmod{q}{p}{n}$, and then, obviously, 
$\xi^{(q)}=\xi^{(p)}$ for all $\xi\in\cA$). We have to prove that
\qtnl{120315a}
c_{X^{(p)}Y^{(p)}}^{Z^{(p)}}=c_{X^{}Y^{}}^{Z^{}}
\eqtn
for all $X,Y,Z\in\cla(G)$, where the numbers in the both sides are the structure constants of  the S-ring $\cZ(\mZ G)$.  
Since this ring is commutative and $p$ is prime, formula~\eqref{120315c} implies that 
$$
\und{X}^{(p)}\,\und{Y}^{(p)}\equiv\und{X}^p\und{Y}^p=
(\und{X}\,\und{Y})^p=(\sum_Zc_{XY}^Z\und{Z})^p\equiv
\sum_Zc_{XY}^Z\und{Z}^{(p)}\ (\text{\rm mod}\hspace{2pt}p).
$$
Thus the relation \eqref{120315a} is true modulo~$p$. Since $p>2n$, we are done.\bull

There is an alternative way to prove the second part of Theorem 1.2. It is related to the action of the group $\mZ_n^*$ on the 
set $\Irr(\cA)$ of all irreducible $\mC$-characters of the S-ring~$\cA$, where as before, we can assume that $\cA=\cZ(\mZ G)$.
Let $\varepsilon$ be an $n$-th primitive complex root of unity. 
Then each $m\in\mZ_n^*$ determines an automorphism~$\tau_m$ of the cyclotomic field $\mQ(\varepsilon)$, 
which sends $\varepsilon$ to $\varepsilon^m$. It follows that  for any $\chi\in\Irr(G)$, the function 
$\chi^{\tau_m}(g):=(\chi(g))^{\tau_m}$, $g\in G$, is also an irreducible character of $G$ and
$$
\chi^{\tau_m}(g) = \chi(g^m)
$$ 
 (see \cite[Proposition~3.16]{Hu}). The primitive idempotents of $\cA$ coincide with the central primitive idempotents of the group algebra $\mQ(\varepsilon)[G]$ which, in turn, are in a one-to-one correspondence with the irreducible characters of $G$. More precisely, if
$e_\chi$ is the idempotent corresponding to $\chi\in \Irr(G)$, then
$$
e_\chi = \frac{1}{|G|}\sum\chi(g)g^{-1}.
$$
A direct computation shows that $\sigma_m(e_\chi) = e_{\chi^{\tau_m}}$. Thus, $\sigma_m$ permutes the primitive idempotents of $\cA$.
This implies that $\sigma_m$ is an automorphism of $\cA$, as required. We note that the above formula shows that there is a natural 
one-to-one correspondence between the $\Irr(G)$ and $\Irr(\cA)$. More precisely, 
\qtnl{170415a}
\Irr(\cA)=\{\frac{1}{\chi(1)}\chi |_{\cA}:\ \chi\in\Irr(G)\}.
\eqtn

Given a set $X\subseteq G$, denote by $\tr(X)$ the union of the sets~$X^{(m)}$, where $m$ runs over the integers coprime to~$n=|G|$;
it is called the {\it trace} of $X$. 
Let  $\cA$ be a central S-ring over $G$. Then from Theorem~\ref{100315v}, it follows that $\und{\tr(X)}\in\cA$ for all
$X\in\cS(\cA)$. Therefore, 
$$
\tr(\cA)=\Span\{\und{\tr(X)}:\ X\in\cS(\cA)\}
$$ 
is a submodule of $\cA$. It is easily seen that it consists of all fixed points of the natural action of the group 
$\{\sigma_m:\ (m,n)=1\}$ on $\cA$. Thus, $\tr(\cA)$ is  an S-ring, which is obviously central; it is called the {\it rational closure} 
of the S-ring~$\cA$. It should be noted that our definitions agreed with the relevant definitions in the abelian case.
The following statement immediately follows from the fact that the $\tr(H)=H$ for any  group $H\le G$.

\prpstn{120315d}
Let $\cA$ be a central S-ring over $G$. Then $\cA$ is primitive if and only if so is $\tr(\cA)$.\bul
\eprpstn
We say that a central S-ring  is {\it rational} if it coincides with its rational closure, or equivalently, 
if each of its  basic sets is rational. The following statement justified the term ``rational''.

\thrml{120315x}
Let $\cA$ be a central S-ring over a group $G$. Then it is rational if and only if  $\pi(\und{X})\in\mQ$ for
all $\pi\in\Irr(\cA)$  and all $X\in\cS(\cA)$.
\ethrm
\proof Let $m$ be an integer  coprime to $n=|G|$. Since any character $\pi\in\Irr(\cA)$ is equal to the restriction to $\cA$ 
of a suitable character $\chi\in\Irr(G)$, from relation~\eqref{170415a} it follows that
\qtnl{210315a}
\pi(\und{X})^{\tau_m}=\pi(\und{X}^{(m)}),\quad X\in\cS(\cA),
\eqtn
where $\tau_m$ is the above defined automorphism of the field $\mQ(\varepsilon)$.
If the S-ring $\cA$ is rational, then the right-hand side of this equality does not depend on the choice of~$m$. So
the number $\pi(\und{X})^\tau$ does not depend on the automorphism $\tau$ of $\mQ(\varepsilon)$. Thus, $\pi(X)\in\mQ$.\medskip

Assume now that  $\pi(\und{X})\in\mQ$ for all $\pi\in\Irr(\cA)$ and $X\in\cS(\cA)$. Then from~\eqref{210315a} it  follows  that
$\pi(\und{X}) = \pi(\und{X}^{(m)})$ for all integers $m$ coprime to $n$ and all characters $\pi\in\Irr(\cA)$. This implies that
$$
\und{X}=\sum_{\pi\in\Irr(\cA)}\pi(\und{X})e_\pi=\sum_{\pi\in\Irr(\cA)}\pi(\und{X}^{(m)})e_\pi=\und{X}^{(m)}),
$$
where $e_\pi$ is the primitive idempotent corresponding to the character~$\pi$. Thus, the S-ring $\cA$ is rational.\bull

\section{Proof of Theorem~\ref {070215a}}\label{190315u}

By the theorem hypothesis, $G$ has a normal Sylow $p$-subgroup $P\cong\mZ_{p^n}$. So by the Schur-Zassenhaus theorem,
$G=PK$ , where $K$ is a Hall $p'$-subgroup of $G$. In what follows, we denote by $H$ the unique subgroup
of $P$ of order~$p$.

\lmml{180315a}
Let $x\in G$ be such that $Hx\not\subset x^G$. Then $x\in C_G(P)$.
\elmm
\proof The element~$x$ acts by conjugation as an automorphism of the cyclic group~$P$. Therefore, 
there exists  an integer~$m$ coprime to $p$ such that $h^x=h^m$ for all $h\in P$.
Rewriting this equality as $x^h =x h^{1-m}$, we obtain 
$$
x^G\supseteq x^P\supseteq P^{(1-m)}x.
$$ 
Since $P^{(1-m)}$ is a subgroup of $P$ and $x^G\not\supseteq xH$, this implies that $P^{(1-m)}=e$. Thus,
$x^h = x$ for all $h\in P$, which means that $x\in C_G(P)$.\bull

Suppose on the contrary that the S-ring $\cA$ is  primitive. Take a nontrivial basic set $X$, which intersects $H$
nontrivially. Then $\grp{X}\ne H$: indeed, otherwise $\grp{X}=H$ by the primitivity of $\cA$ and $n=p$ is a prime
in contrast to the hypothesis. This proves the second 
part of the following relations (the first one follows from the choice of~$X$):
\qtnl{140315c}
X\cap H\ne\varnothing\qaq X\setminus H\ne\varnothing\qaq \grp{X\cap H}\le\rad(X\setminus H).
\eqtn
To prove the third one, set $X_0=\{x\in X:\ xH\not\subset X\}$. Then from Lemma~\ref{180315a} it follows that
\qtnl{130315b}
(X_0)^{(p)}\subseteq P^{(p)}K\subsetneq G.
\eqtn
Moreover, it is easily seen that the sets  $X_0$ and $X\setminus X_0$ are unions of some conjugacy classes of $G$. For
these classes, we can refine Lemma~\ref{100315b} as follows.

\lmml{250315a}
For any class $Y\in\cla(G)$, we have 
$$
\und{Y}^p\equiv\css
\und{Y}^{(p)}\ (\text{\rm mod}\hspace{2pt}p), &\text{if $Y\subseteq C_G(P)$,}\\
\ 0\quad\ (\text{\rm mod}\hspace{2pt}p), &\text{if $Y\not\subseteq C_G(P)$.}\\
\ecss
$$
\elmm
\proof The group $C:=C_G(P)$ is obviously normal in $G$. Therefore, 
$$
Y\subseteq C\qoq Y\cap C=\varnothing. 
$$
Suppose first 
that $Y\cap C=\varnothing$. Since $H\trianglelefteq G$, we have $y\und{H}=\und{H}y$ for all~$y\in Y$.
Denote by $S$ a full system of representatives of the family $\{Hy:\ y\in Y\}$. Then, since $|H|=p$, we have
$$
\und{Y}^p=(\sum_{y\in S}\und{H}y)^p\equiv\und{H}^p\und{S}^p\equiv 0\ (\text{\rm mod}\hspace{2pt}p),
$$
as required. Let now $Y\subseteq C$. Then $Y$ is a normal subset of $C$, i.e. $Y^G=Y$.  Since the group~$C$ is a direct product of $P$ 
and $O_{p'}(C)$, each normal subset of $C$ is the disjoint union of $gY_g$, $g\in P$, where $Y_g$ is a normal subset of $O_{p'}(C)$. 
Now 
\qtnl{180415a}
\underline{Y}^{(p)}=\sum_{g\in P} \underline{gY_g}^{(p)} 
=\sum_{g\in P} g^p\,\underline{Y_g}^{(p)}.
\eqtn
Moreover, since $Y_g$ is contained in the $p'$-subgroup $O_{p'}(C)$,  by Lemma~\ref{100315b} we obtain
\qtnl{180415b}
\underline{Y_g}^{(p)} = \underline{Y_g^{(p)}}\equiv \underline{Y_g}\,^p\ (\text{\rm mod}\hspace{2pt}p).
\eqtn
Thus, from \eqref{180415a} and \eqref{180415b}, it follows that
$$
\underline{Y}^p\equiv\sum_{g\in P} (g\underline{Y_g})^p=\sum_{g\in P} g^p \underline{Y_g}^p\equiv 
\sum_{g\in P} g^p\,\underline{Y_g}^{(p)}=\underline{Y}^{(p)}\ (\text{\rm mod}\hspace{2pt}p)
$$ 
as required.\bull

To complete the proof of the third relation in~\eqref{140315c}, suppose on the contrary that the set~$X_0$ is not empty. 
Then, if $X$ is the union of conjugacy classes $X_i$, $i\in I$, then by Lemma~\ref{250315a}, we have
\qtnl{250315r}
\und{X}^p=(\sum_{i\in I}\und{X_i})^p=\sum_{i\in I}\und{X_i}^p=\sum_{i\in I_0}\und{X_i}^{(p)}\ (\text{\rm mod}\hspace{2pt}p),
\eqtn
where $I_0=\{i\in I:\ X_i\subseteq X_0\}$. 
Moreover, by Lemma~\ref{180315a}, given $x,y\in X_0$, the equality $x^p=y^p$ holds  if and only if $y\in Hx$. Since
also 
$$
1\le |Hx\cap X_0|\le p-1
$$ 
for all $x\in X_0$, the coefficient at $x^p\in G$ in the right-hand sum of~\eqref{250315r} is 
between $1$ and $p-1$. Thus, 
$$
\xi:=f[\und{X}^p]
$$ 
is a non-zero element of the S-ring~$\cA$, where $f$ is the function used
in~\eqref{110315a}. By the Schur-Wielandt principle, this implies that the support~$Y$ of the element~$\xi$ is an $\cA$-set.
Therefore, $\grp{Y}$ is an $\cA$-subgroup of~$G$. This subgroup is proper: $\grp{Y}\ne G$
by~\eqref{130315b} and $\grp{Y}\ne e$, because $X_0\ne\varnothing$. 
But this contradicts the primitivity of the S-ring~$\cA$.\medskip

Thus, all the relations in~\eqref{140315c} are true. To complete the proof,  we make use of the following theorem  on separating subgroup 
proved in~\cite{EP}.

\thrml{t100703}
Let $\cA$ be an   S-ring over a group $G$. Suppose that $X\in\cS(\cA)$ and $H\le G$ satisfy relations~\eqref{140315c}.
Then $X=\grp{X}\setminus\rad(X)$ and $\rad(X)\le H\le\grp{X}$.\bul
\ethrm

Now, since $\rad(X)$ and $\grp{X}$ are $\cA$-groups,  the primitivity assumption implies that
$\rad(X)=e$ and $\grp{X}=G$. By Theorem~\ref{t100703}, this implies that $X=G\setminus e$. This means that 
$\rk(\cA)=2$, i.e., the S-ring $\cA$ is trivial. Contradiction.

\section{Generalized B-groups}\label{190315v}

\sbsnt{Proof of Theorem~\ref{100315a}.}
Let $G$ be a Camina group. Then it has a normal subgroup $H$ such that $(G,H)$ is a Camina pair.
Let $\cA$ be a proper central primitive S-ring over  $G$. Take a set $X\in\cS(\cA)$ that contains
a nonidentity element of~$H$. It follows from the primitivity of $\cA$ that
\qtnl{250315f}
\rad(X)=e\qaq \grp{X}=G.
\eqtn
In particular, the first two relations in~\eqref{140315c} hold. Next, the set
$X$ is a union of some conjugacy classes of~$G$ as the S-ring $\cA$ is central. By the definition
of a Camina pair, we have
$$
xH=Hx\subseteq X\setminus H
$$
for all $x\in X\setminus H$. This proves the third relation in~\eqref{140315c}. Thus, $X=\grp{X}\setminus\rad(X)$
by Theorem~\ref{t100703}. By~\eqref{250315f}, this implies that $X=G\setminus e$ and hence $\rk(\cA)=2$. The latter 
means that the S-ring~$\cA$ is not proper. Contradiction.\bull

\sbsnt{A generalized B-group, which is not a B-group.}\label{190315a}
Let $p>3$ be a prime congruent to $3$ modulo $4$, and let $G$ be the extraspecial group of order $p^3$ and exponent $p$.
Then there exists  a skew Hadamard difference set $X$ in the group $G$; see~\cite{Fe}. This exactly means that $Y:=X^{-1}$
is equal to $G^\#\setminus X$ and
$$
\und{X}\und{Y}=|X|e + \frac{|X|-1}{2}(\und{X}+\und{Y}).
$$
Therefore, the module $\cA=\Span\{e,\und{X},\und{Y}\}$ is a subring of $\mZ G$ that satisfies the conditions
(S1), (S2), and (S3) with $\cS=\{e,X,Y\}$. Thus, $\cA$ is an S-ring of rank~$3$ over $G$. This S-ring is, obviously,
primitive. Since it is also proper, $G$ is not a B-group. On the other hand, it is a generalized B-group by Corollary~\ref{100315u}.

\sbsnt{Simple groups.} 
According to~\cite{FS}, a group $G$ is said to be {\it rational} if the number $\chi(g)$ is rational for all 
$\chi\in\Irr(G)$ and all $g\in G$. Finite simple rational groups were characterized in Corollary~B1 of that paper
as follows: a noncycic simple group $G$ is rational if and only if $G\cong\SP(6,2)$ or $\ORT^+(8,2)'$.\medskip

{\bf Proof of Theorem~\ref{070215b}.} Let $G$ be a finite simple group other than $\SP(6,2)$ or $\ORT^+(8,2)'$. Without loss of 
generality, we can assume that $G$ is not cyclic. Then by the above characterization of rational groups, $G$ is not rational and has
two elements $x$ and $y$ of distinct orders. Then the orders of $x^m$ and $y^m$ are also distinct for all integers $m$ coprime to~$|G|$. 
This implies that the order of any element of $\tr(x^G)$ does not equal the order of any element of $\tr(y^G)$. So,
\qtnl{250315s}
\tr(x^G)\ne \tr(y^G).
\eqtn
Therefore, the rational closure $\tr(\cA)$ of the S-ring $\cA=\cZ(\mZ G)$ is
of rank at least~$3$. On the other hand, $\tr(\cA)\ne\cA$, 
for otherwise the irreducible characters of $\cA$ are rational valued (Theorem~\ref{120315x}) and then $G$ is a rational group.
Thus, $\tr(\cA)$ is a proper central S-ring. It is primitive because so is $\cA$ (Proposition~\ref{120315d}). Therefore $G$
can not be a generalized B-group.\bul

\sbsnt{AS-free groups.}
According to~\cite{ABC}, a  transitive permutation group is called {\it AS-free} if it preserves no nontrivial 
symmetric association scheme. From Theorem~17 of that paper, it follows that given  a nonabelian  simple group $G$,
the permutation group on $G$ defined by
$$
K=\grp{G_{right},\aut(G),\sigma},
$$ 
is AS-free, where $G_{right}$ is the group of all
 right translations of $G$ and $\sigma$ is a permutation of $G$ that takes $g$ to $g^{-1}$, $g\in G$. It is easily
seen that the orbits of the stabilizer of $e$ in $K$ are the basic sets of a central S-ring $\cA$ over $G$. Thus, using the
above result one can get another proof that $G$ is not a generalized B-group whenever $\cA\ne \cZ(\mZ G)$. However,
in general, the latter inequality is not true, e.g., for the group $\SP(6,2)$.

\sbsnt{Miscellaneous.}
Let $G$ be a finite group having a relatively prime conjugacy class (examples of such groups can be found, e.g., in~\cite{DMN}). 
Denote by $\cX$ the association scheme of the permutation group $G\Inn(G)\le\sym(G)$ (see also, \cite[Theorem~7.2]{BI}).
Then one can see that $\cla(G)$ forms a relatively prime equitable partition for $\cX$ in the sense of~\cite{HKK}.
By Theorem~3.1 of that paper, any primitive fusion of the scheme $\cX$ must have rank~$2$. So, using the correspondence
between the Cayley schemes and S-rings over $G$, one can show that $G$ is a generalized B-group.\medskip

Let $G=G_1\times G_2$ where $G_1$ and $G_2$ are groups of the same order $n>1$. Then $G$ is not a generalized B-group.
Indeed, set $X_0=\{(e_1, e_2)\}$ where $e_i$ is the identity of $G_i$, and
$$
X_1=e_1\times (G_2)^\#\,\cup (G_1)^\#\times e_2.
$$
Denote by $\cA$ the span of the set $\{\und{X_i}:\ i=0,1,2\}$ where $X_2$ is the complement to $X_0\cup X_1$ in $G$.
Then  $\cA$ is, obviously, a central S-ring of rank~$3$ over $G$. Since it is also primitive, we are done.\medskip

\centerline{\bf Acknowledgment.}
\medskip

The paper was started during the visit of the third author to the Central China Normal University, Wuhan, China. He would like
to thank the faculty members of the School of Mathematics and Statistics for their hospitality.

\end{document}